\documentclass[12pt, a4paper, reqno, oneside]{amsart}
\usepackage{amsmath,amsthm,amssymb,mathrsfs}
\usepackage{color,hyperref,graphicx,comment}
\usepackage[utf8]{inputenc}
\usepackage[all]{xy}
\usepackage{mathtools}
\usepackage{enumitem}
\usepackage{faktor}
\usepackage{tikz-cd}
\usepackage[all]{xy} %xymatrix
\usepackage{upgreek}
\usepackage{bbold} %\mathbb{1}
\usepackage{stmaryrd} %\owedge
\usepackage{indentfirst}
\usepackage{empheq}
\usepackage{lipsum}
\usepackage{natbib}
\setcitestyle{numbers,square}
\usepackage{geometry}
\theoremstyle{thmstyleone}
\newtheorem{theorem}{Theorem}

\theoremstyle{thmstyletwo}

\newtheorem{defn}{Definition}

\newtheorem{cor}{Corollary}
\numberwithin{equation}{section}
\numberwithin{defn}{section}
\numberwithin{lemma}{section}
\numberwithin{cor}{section}
\numberwithin{example}{section}

\def\R{\mathbb{R}}

\def\S{\mathbb{S}}

\begin{document}
\title{Translation Surfaces in Lie Groups with Constant Gaussian Curvature}
\author{Xu Han and Zhonghua Hou}
\address{Dalian University of Technology, Dalian City, China, 116024}
\email{xuhan2016@mail.dlut.edu.cn and zhhou@dlut.edu.cn}

\begin{abstract}
\noindent
Let $G$ be a $n$-dimensional Lie group with a bi-invariant Riemannian metric. We prove that if a surface of constant Gaussian curvature in $G$ can be expressed as the product of two curves, then it must be flat. In particular, we can essentially characterize all such surfaces locally in 3-dimensional case.
\end{abstract}

\maketitle

\section{Introduction}\label{sec1}

In \cite{Hasa2021}, Hasanis and Lopez define a surface $M$ in $\R^3$ as a translation surface if it can be locally written as the sum of two space curves: $f(s,t)=\alpha(s)+\beta(t)$. They prove that any regular translation surface of constant Gaussian curvature must be flat. The historical motivation and the research background of the translation surfaces are well introduced in their paper and references therein.

In this paper, we generalize the concept to the following

\begin{defn}
A surface $M$ in an n-dimensional Lie group $G$ is called a translation surface if it can be locally written as the product of two curves: $f(s,t)=a(s)\cdot b(t)$, where the product is the multiplication of $G$.
\end{defn}
\noindent
And we prove a similar result

\begin{theorem}
Let $G$ be an n-dimensional ($n\geq 2$) Lie group with a bi-invariant metric, and $S$ be a regular translation surface in $G$ with constant Gaussian curvature, then the surface must be flat.
\end{theorem}
\noindent
\textbf {Remark}. Particularly, $\R^n$ can be regarded as a Lie group with sum being the group multiplication. And when $n=2$, the above theorem is trivial, since the only simply-connected two dimensional Lie group with bi-invariant metric is the Euclidean space $\R^2$. To prove Theorem 1, we apply the method of symmetry group theory, without recourse to the complex calculation of the resultant in \cite{Hasa2021}.

Since it is well known (see \cite{Mil1976} for example) that the universal covering groups of the connected 3-dimensional Lie groups with bi-invariant metric can only be $\R^3$ and $\S^3$, others can be obtained by the quotients of normal discrete subgroups, i.e. $S^1\times\R^2$, $S^1\times S^1\times\R^1$, $S^1\times S^1\times S^1$, $SO(3)$. To characterize the local properties of translation surfaces, it suffices to study $\R^3$ and $\S^3$. According to \cite{Hasa2021} , we can easily derive that besides cylindrical surfaces, locally there is only one extra possibility of translation surfaces with constant Gaussian curvature, that is a part of a flat tori in $S^1\times S^1\times\R^1$ and $S^1\times S^1\times S^1$. As for $\S^3$, we have the following characterization due to Spivak's work \cite{Spiv1975}

\begin{cor}
Let $M$ be a translation surface of constant Gaussian curvature in $\S^3$, for any generic point $p\in M$, there exists a neighbourhood $U\subset M$ of $p$ such that $U$ can be written as the product of two asymptotic curves, whose torsion are 1 and -1 respectively. Moreover, if the translation surface is further simply-connected and complete, then such an expression is global.
\end{cor}
\noindent
By generic point we mean a point where the curvatures of both asymptotic curves are not zero. Further details of Spivak's work are recalled in Section 2.

\section{Preliminaries}
For the convenience, we include in this section the necessary facts about symmetry group theory and flat surfaces in the unit 3-sphere, and suggest the reader refer to \cite{Olv1993} and \cite{Spiv1975} for the detailed accounts. \par

\subsection{Symmetry group}

\begin{defn}
Let $M$ be a smooth manifold and $G$ be a Lie group or a local Lie group. Suppose that $\mathcal{U}$ is an open subset of $G\times M$, such that
\begin{equation}
\{e\}\times M\subset \mathcal{U}.
\end{equation}
A local group of transformations acting on $\mathcal{U}$ is defined as a smooth map $\Omega:\mathcal{U}\rightarrow M$ with the following properties:\par
\noindent
(a) If $(h,x)\in\mathcal{U}$, $(g,\Omega(h,x))\in\mathcal{U}$ and $(g\cdot h,x)\in \mathcal{U}$, then
\begin{equation}
\Omega(g,\Omega(h,x))=\Omega(g\cdot h,x).
\end{equation}
(b) For all $x\in \mathcal{U}$,
\begin{equation}
\Omega(e,x)=x.
\end{equation}
(c) If $(g,x)\in \mathcal{U}$, then $(g^{-1}, \Omega(g,x))\in\mathcal{U}$ and
\begin{equation}
\Omega(g^{-1},\Omega(g,x))=x.
\end{equation}
\end{defn}

\begin{defn}
Let $G$ be a local group of transformations defined on some open subset $U\subset\{(x,u)\in\mathbb{R}^n\times\mathbb{R}^q\}$, which is the space of independent and dependent variables of a PDE system. We call $G$ is a symmetry group of the system, if $u=f(x)$ is a solution to the system, then whenever $g\cdot f$ is defined for $g\in G$, we have $u=g\cdot f(x)$ is also a solution.
\end{defn}
\noindent
\textbf{Remark.} We only consider the case $n=2$ and $q=1$ in the paper.

In order to calculate the symmetry group of a PDE, we first extend the usual space of variables to the jet space to include partial derivatives, which transforms the original system into a system of algebraic equations. We now give the definition of prolongation.
\begin{defn}
The $m$-th prolongation of a smooth function $u=f(x)$, which is denoted by $u^{(m)}=pr^{(m)}f(x)$, is an vector-valued function $\{u_{J}\}$ for all multi-indices $J=(j_1,...,j_k)$, and $0\leqslant k\leqslant m$ defined by the following equations
\begin{equation}
u_{J}:=\frac{\partial^{k} f(x)}{\partial x^{j_1}\cdot\cdot\cdot\partial x^{j_k}},
\end{equation}
particularly, we define $u_{J}=u$ when $k=0$.
\end{defn}

\begin{defn}
Given a local group of transformations $G$ acting on an open subset $U$ of the variable space, and for any smooth function $u=f(x)$ such that $u_0=f(x_0)$ and $(\tilde x_0, \tilde u_0)=g\cdot(x_0,u_0)$ is defined for $g\in G$, we define the $m$-th prolongation of $g$ at the point $(x_0,u_0)$ by
\begin{equation}
pr^{(m)}g\cdot(x_0,u_0^{(m)})=(\tilde x_0, \tilde {u}_0^{(m)}),
\end{equation}
where
\begin{equation}
\tilde {u}_0^{(m)}:=pr^{(m)}(g\cdot f)(\tilde x_0).
\end{equation}
\end{defn}

\begin{defn}
Suppose $v$ is a vector field on $U\subset\mathbb{R}^{n+1}$, with the corresponding (local) one-parameter transformations $\Omega_{\varepsilon}$. Then $m$-th prolongation of $v$ at the point $(x,u^{(m)})$ , which is denoted by $pr^{(m)}v(x,u^{(m)})$, is defined by the following equation
\begin{equation}
pr^{(m)}v(x,u^{(m)})=\frac{d}{d\varepsilon}\Big |_{\varepsilon=0}pr^{(m)}\Omega_{\varepsilon}(x,u^{(m)})
\end{equation}
\end{defn}

\begin{defn}
\label{def max rank}
A system of differential equations
\begin{equation}
\Phi^{\alpha}(x,u^{(m)})=0,\quad \alpha=1,...,l
\end{equation}
is of maximal rank if the Jacobian matrix of $\Phi:=\{\Phi^{\alpha}\}$
\begin{equation}
J_{\Phi}(x,u^{(m)})=\{(\frac{\partial \Phi^{\alpha}}{\partial x^{i}}),(\frac{\partial \Phi^{\alpha}}{\partial u_{J}})\},
\end{equation}
is of rank $l$ on $\mathcal{S}:=\{(x,u^{(m)}):\Phi(x,u^{(m)})=0\}$.
\end{defn}
\noindent
\textbf{Remark.} We only consider the case $l=1$ in the paper.

\begin{defn}
\label{def loca solv}
A system of differential equations $\Phi(x,u^{(m)})=0$ is said to be locally solvable at a point $(x_{0}, u_{0}^{(m)})\in \mathcal{S}$, if there is a solution $u=f(x)$ of the system defined in a neighbourhood of $x_{0}$, such that $u_{0}^{(m)}=pr^{(m)}f(x_{0})$

\end{defn}

For the partial differential equations in Kovalevskaya form, the Cauchy-Kovalevskaya theorem provides us with a criterion to check the local solvability.

\begin{theorem}
\label{cauchy-kova}
An analytic PDE in the following Kovalevskaya form is locally solvable:
\begin{equation}
\frac{\partial^{m}u}{\partial t^m}=F(x,t,\widetilde{u^{(m)}}),
\end{equation}
where $(x,t)=(x^1,...,x^{n-1},t)$ are the independent variables, and $\widetilde{u^{(m)}}$ denotes all partial derivatives of $u$ with respect to $x$ and $t$ up to order $n$ except $\frac{\partial^{m}u}{\partial t^m}$.
\end{theorem}

\begin{defn}
A system of differential equations is called non-degenerate, if it is both locally solvable and of maximal rank at every point $(x_{0}, u_{0}^{(m)})\in \mathcal{S}$.

\end{defn}

We have the following necessary and sufficient condition for a group to be the symmetry group.

\begin{theorem}
\label{equi cond}
Let
\begin{equation}
\Phi^{\alpha}(x,u^{(m)})=0,\quad \alpha=1,...,l
\end{equation}
be a non-degenerate system of differential equations defined on $U\subset \mathbb{R}^{n+1}$. If $G$ is a connected local group of transformations acting on $U$, then $G$ is a symmetry group of the system if and only if for every infinitesimal generator $v$ of $G$, the following are satisfied on $\mathcal{S}$,
\begin{equation}
\label{theorem 2 cond}
pr^{(m)}v[\Phi^{\alpha}(x,u^{(m)})]=0,
\end{equation}
for $\alpha=1,...,l$.
\end{theorem}

We also need a prolongation formula.
\begin{defn}
The $i$-th total derivative of a given function $P(x,u^{(m)})$ is defined by the following equation
\begin{equation}
D_iP=\frac{\partial P}{\partial x^i}+\sum_J u_{J,i}\frac{\partial P}{\partial u_J},
\end{equation}
where $J=(j_1,...,j_k)$, and
\begin{equation}
u_{J,i}=\frac{\partial u_J}{\partial x^i}=\frac{\partial^{k+1} u}{\partial x^i \partial x^{j_1}\cdot\cdot\cdot\partial x^{j_k}}.
\end{equation}
The $J$-th total derivative is defined by
\begin{equation}
D_J=D_{j_1}D_{j_2}\cdot\cdot\cdot D_{j_k}.
\end{equation}
\end{defn}

\begin{theorem}
\label{formula}
Let
\begin{equation}
v=\sum_{i=1}^{n}\xi^{i}(x,u)\frac{\partial}{\partial x^{i}}+\eta(x,u)\frac{\partial}{\partial u}
\end{equation}
be a vector field defined on the variable space. Suppose that the $m$-th prolongation of $v$ is
\begin{equation}
pr^{(m)}v=v+\sum_{J}\phi^{J}(x,u^{(m)})\frac{\partial}{\partial u_{J}}.
\end{equation}
Then the coefficient functions are given by
\begin{equation}
\phi^{J}(x,u^{(m)})=D_{J}\bigg(\eta-\sum_{i=1}^{n}\xi^{i}u_{i}\bigg)+\sum_{i=1}^{n}\xi^{i}u_{J,i}.
\end{equation}
\end{theorem}

\subsection{Flat surfaces in the 3-sphere}

Given a smooth flat surface $M$ in the unit 3-sphere $\S^3$, the Gauss-Kronecker of $M$ is -1 by Gauss equation, so for any point $p\in M$, there are exactly two local asymptotic curves through $p$, which form the so-called Chebyshev net by choosing the arc-length parameters for both curves.
\begin{defn}
A point $p$ in $M$ is called a generic point, if the curvatures of the two asymptotic curves are not zero at $p$.
\end{defn}
\noindent
Spivak in \cite{Spiv1975} obtains that

\begin{theorem}
A flat surface $M\subset\S^3$ is locally a translation surface at generic points. Moreover, the translation surface can be written as the product of two asymptotic curves, whose torsion are 1 and -1 respectively.
\end{theorem}

\begin{theorem}
If a flat surface $M\subset\S^3$ is complete and simply-connected, then Chebyshev net is globally defined.
\end{theorem}

Conversely, if we have two curves $\alpha(s)$ and $\beta(t)$ with non-zero curvatures and torsion $\tau(s)=1$ and $\tau(t)=-1$ respectively. Without loss of generality, we assume that $\alpha(0)=\beta(0)=1\in\S^3$ (the unit element of $\S^3$) and they are not tangent at $1$, then the local translation surface $f(s,t)={\alpha(s)\cdot\beta(t)}$ is flat. Moreover, if one or both of the curves are geodesic, we have the same conclusion.\par

\section{Proof of Theorem 1}

Let $(G,\langle\cdot, \cdot\rangle)$ be $n$-dimensional ($n\geq2$) Lie group with a bi-invariant metric $\langle\cdot, \cdot\rangle$, and $a(s)$ and $b(t)$ be  two smooth curves in $G$ with arc-length parameters whose tangent vector fields are $T_{a}$ and $T_{b}$ respectively. Suppose $f(s,t)=a(s)\cdot b(t)$ is a regular translation surface, then
\begin{equation}
f_s(s,t)=(R_{b(t)})_{*a(s)}(T_a(s))
\end{equation}
and
\begin{equation}
f_t(s,t)=(L_{a(t)})_{*b(s)}(T_b(s)).
\end{equation}

Since the metric of $G$ is bi-invariant, the induced metric $g$ is
\begin{equation}
g=ds^2+2\langle f_s,f_t\rangle dsdt+dt^2.
\end{equation}
By Cauchy inequality, we have
\begin{equation}
\langle f_s,f_t\rangle\leq 1
\end{equation}
so we can set
\begin{equation}
\cos{u}:=\langle f_s,f_t\rangle, \quad 0<u<\pi.
\end{equation}
Thus the Gaussian curvature is determined by the metric, and computed as
\begin{equation}
\label{curvature euqation 1}
K=-\frac{u_{st}}{\sin{u}},
\end{equation}
where $u_{st}$ is the second order partial derivatives.

From now on, we assume that $K$ is a constant, and set
\begin{equation}
\label{curvature euqation 2}
\Phi(s,t,u^{(2)}):=u_{st}+K\sin{u}.
\end{equation}
We calculate the symmetry group of $\Phi=0$, which is a routine and easy to carry out in this case. Applying Theorem \ref{formula} to
\begin{equation}
v:=\tau(s,t,u)\frac{\partial}{\partial t}+\xi(s,t,u)\frac{\partial}{\partial s}+\eta(s,t,u)\frac{\partial}{\partial u},
\end{equation}
we obtain the prolongation of $v$.
\begin{equation}
pr^{(2)}v=v+\phi^s\frac{\partial}{\partial u_s}+\phi^t\frac{\partial}{\partial u_t}+\phi^{ss}\frac{\partial}{\partial u_{ss}}+\phi^{st}\frac{\partial}{\partial u_{st}}+\phi^{tt}\frac{\partial}{\partial u_{tt}},
\end{equation}
where the only needed coefficient is
\begin{equation}
\begin{aligned}
\phi^{st}=&\eta_{ts}+\eta_{tu}u_s+\eta_{us}u_t-\xi_{ts}u_s-\tau_{ts}u_t+\eta_{uu}u_su_t-\xi_{tu}u_s^2-\xi_{us}u_tu_s\\
          &-\tau_{tu}u_su_t+\eta_uu_{ts}-\xi_tu_{ss}-\xi_su_{st}-\tau_tu_{ts}-\xi_{uu}u_tu_s^2-\tau_{us}u_t^2\\
          &-2\xi_uu_{ts}u_s-\xi_uu_tu_{ss}-\tau_{uu}u_su_t^2-2\tau_uu_tu_{ts}-\tau_su_{tt}-\tau_uu_su_{tt}.
\end{aligned}
\end{equation}
By Theorem \ref{equi cond}, we apply $pr^{(2)}v$ to (\ref{curvature euqation 2}) to obtain
\begin{equation}
\label{prv=0}
\phi^{st}+K\eta\cos u=0,
\end{equation}
\noindent
which is satisfied whenever $\Phi=0$. This forms a system of algebraic equations where the dependence of monomials is restricted to the condition $\Phi(s,t,u^{(2)})=0$. \par
Straightforward calculation yields
\begin{equation}
\begin{aligned}
0=&-\xi_{uu} u_s^2 u_t-\xi_{tu} u_s^2-\xi_t u_{ss}-\xi_u u_{ss} u_t\\
  &-\tau_{uu} u_s u_t^2-\tau_u u_s u_{tt}-\tau_s u_{tt}-\tau_{us} u_t^2\\
  &+u_su_t\left(\eta_{uu}-\xi_{us}-\tau_{tu}\right)\\
  &+u_s \left(2 K \xi_u \sin u+\eta_{tu}-\xi_{ts}\right)\\
  &+u_t \left(2 K \tau_u\sin u+\eta_{us}-\tau_{ts}\right)\\
  &+\eta_{st}+K \eta_u \sin u-K \eta \cos u+K \xi_s \sin u+K \tau_t \sin u.
\end{aligned}
\end{equation}
Taking the coefficients of various independent monomials in partial derivatives of $u$, we derive the determining equations for the symmetry group

\begin{subequations}
\label{determing equation}
\begin{empheq}[left=\empheqlbrace]{align}
&\xi=\xi(s)\\
&\tau=\tau(t)\\
&\eta_{uu}-\xi_{us}-\tau_{tu}=0\\
&2 K \xi_u \sin u+\eta_{tu}-\xi_{ts}=0\\
&2 K \tau_u\sin(u)+\eta_{us}-\tau_{ts}=0\\
&\eta_{st}-K \eta_u \sin u+K \eta \cos u+K \xi_s \sin u+K \tau_t \sin u=0,
\end{empheq}
\end{subequations}
We can solve this system of equations directly and claim that if $K\neq0$, then $\eta=0$. In fact, by ($\ref{determing equation}a$), ($\ref{determing equation}b$) and ($\ref{determing equation}c$), we have $\eta=\alpha(s,t)u+\beta(s,t)$, for some functions $\alpha$ and $\beta$. Combining ($\ref{determing equation}d$) and ($\ref{determing equation}e$), we see that $\alpha$ is a constant, i.e. $\eta=\alpha u+\beta(s,t)$. Inserting this into ($\ref{determing equation}f$) yields
\begin{equation}
\label{equ f}
\beta_{st}-K \alpha \sin u+K (\alpha u+\beta(s,t))\cos u+K \xi_s \sin u+K \tau_t \sin u=0.
\end{equation}
Since $K\neq0$, taking the derivative of ($\ref{equ f}$) with respect to $u$ yields
\begin{equation}
\label{u deriv}
(\xi_s+\tau_t) \cos u-(\alpha u+\beta(s,t))\sin u=0.
\end{equation}
\noindent
Again taking the derivatives of ($\ref{u deriv}$) with respect to $s$ and $t$ respectively, i.e.
\begin{equation}
\begin{cases}
\xi_{ss}\cos u-\beta_s\sin u=0\\
\tau_{tt}\cos u-\beta_{t}\sin u=0,
\end{cases}
\end{equation}
we obtain the following also due to $K\neq0$, which implies that $u$ is not constant.
\begin{equation}
\label{const equ}
\begin{cases}
\xi=As+B\\
\tau=Ct+D\\
\beta=E,
\end{cases}
\end{equation}
where $A$, $B$, $C$, $D$ and $E$ are constants. Substitute (\ref{const equ}) into (\ref{equ f}), we see $A+C=0$ and $\alpha=\beta=0$. Thus we prove the claim.\par

Whether $K=0$ or not, in order to apply Theorem  \ref{equi cond} to verify the symmetry group of the original PDE, we have to check the non-degeneracy condition. The maximal rank condition is trivially satisfied, and the local solvability condition can also be easily checked by Theorem \ref{cauchy-kova} considering the following change of variables, which transforms (\ref{curvature euqation 1}) into Kovalevskaya form.
\begin{equation}
\begin{cases}
x=t+s\\
y=t-s.
\end{cases}
\end{equation}

Therefore, we obtain that when $K\neq0$, the infinitesimal symmetries of (\ref{curvature euqation 1}) have the following form
\begin{equation}
\label{v form}
v=(-At+D)\frac{\partial}{\partial t}+(As+B)\frac{\partial}{\partial s}
\end{equation}

However, the $s$-curve or $t$-curve induces a local one-parameter family of diffeomorphisms, which preserves (\ref{curvature euqation 1}). By definition it belongs to the symmetry group. In fact, if we assume that for some short interval $I$ containing 0, and $\epsilon\in I$, then the local one-parameter family of diffeomorphisms induced by the $s$-curve acting on the initial local surface $f(s,t)=a(s)\cdot b(t)$ near a point $p\in M$ is
\begin{equation}
f_{\epsilon}:=a(s+\epsilon)\cdot b(t).
\end{equation}
So
\begin{equation}
\begin{aligned}
\langle \frac{\partial f_{\epsilon}}{\partial s},\frac{\partial f_{\epsilon}}{\partial t}\rangle&=\cos u_{\epsilon}\\
\langle \frac{\partial f_{\epsilon}}{\partial s},\frac{\partial f_{\epsilon}}{\partial s}\rangle&=1\\
\langle \frac{\partial f_{\epsilon}}{\partial t},\frac{\partial f_{\epsilon}}{\partial t}\rangle&=1.
\end{aligned}
\end{equation}
where $u_{\epsilon}$ is the induced transformation of $u$. Taking derivative with respect to $\epsilon$ and by (\ref{v form}), we have when $K\neq0$,
\begin{equation}
\frac{du_{\epsilon}}{d\epsilon}{\Big|}_{\epsilon=0}=0.
\end{equation}
Thus along $s$-curve $u$ should be constant, since the point $p$ is arbitrary, it implies that the Gaussian curvature is zero, which is a contradiction.

In conclusion, we prove that the constant Gaussian curvature can only be zero, and the infinitesimal symmetry have the following form

\begin{equation}
v=\tau(t)\frac{\partial}{\partial t}+\xi(s)\frac{\partial}{\partial s}+(\alpha u+\gamma(s)+\delta(t))\frac{\partial}{\partial u}.
\end{equation}

\bibliographystyle{abbrv}

\begin{thebibliography}{10}


\bibitem{Hasa2021} Thomas Hasanis and Rafael Lopez, Translation surfaces in Euclidean space with constant Gaussian curvature, \textit{Comm. Anal. Geom.}, 29(2021), No. 6, 1415-1447.
\bibitem{Mil1976} John Milnor, Curvatures of left invariant metrics on Lie groups,  adv. math., 21(1976), 293-329.
\bibitem{Olv1993} Peter Olver, Applications of Lie groups to the differential equations, Springer, 2nd Edition, 1993.
\bibitem{Spiv1975} Michael Spivak, A Comprehensive introduction to differential geometry, Publish or Perish, Boston, Vol. IV, 1975.

\end{thebibliography}

\end{document}